
\def\title#1{{\titlefont\noindent #1\bigskip}}

\def\author#1{{\largefont\noindent #1}\medskip}

\def\beginlinemode{\endmode
 \begingroup\obeylines\def\endmode{\par\endgroup}}
\let\endmode=\par

\newbox\theaddress
\def\address{\smallskip\beginlinemode\parindent 0in\getaddress}
{\obeylines
\gdef\getaddress #1 
 #2
 {#1\gdef\addressee{#2}%
   \global\setbox\theaddress=\vbox\bgroup\raggedright%
    \everypar{\hangindent2em}#2
   \def\endaddress{\egroup\endgroup \copy\theaddress \medskip}}}

\def\thanks#1{\footnote{}{\eightpoint #1}}

\long\def\Abstract#1{{\eightpoint\narrower\vskip\baselineskip\noindent
#1\smallskip}}

\def\skipfirstword#1 {}

\def\ir#1{\csname #1\endcsname}

\newdimen\currentht
\newbox\droppedletter
\newdimen\droppedletterwdth
\newdimen\drophtinpts
\newdimen\dropindent

\def\irrnSection#1#2{\edef\tttempcs{\ir{#2}}
\vskip6pt\penalty-3000
{\bf\noindent \expandafter\skipfirstword\tttempcs. #1}
\vskip6pt}

\def\irSubsection#1#2{\edef\tttempcs{\ir{#2}}
\vskip\baselineskip\penalty-3000
{\bf\noindent \expandafter\skipfirstword\tttempcs. #1}
\vskip6pt}

\def\irSubsubsection#1#2{\edef\tttempcs{\ir{#2}}
\vskip\baselineskip\penalty-3000
{\bf\noindent \expandafter\skipfirstword\tttempcs. #1}
\vskip6pt}

\def\References{\vbox to.25in{\vfil}\noindent{}{\bf References}
\vskip6pt\par}

\def\References{\vskip6pt\noindent{}{\bf References}
\vskip6pt\par}

\def\baselinebreak{\par \ifdim\lastskip<6pt
         \removelastskip\penalty-200\vskip6pt\fi}

\long\def\prclm#1#2#3{\baselinebreak
\noindent{\bf \csname #2\endcsname}:\enspace{\sl #3\par}\baselinebreak}

\def\Prf{\noindent{\bf Proof}: }

\def\rem#1#2{\baselinebreak\noindent{\bf \csname #2\endcsname}: }

\def\qed{{$\diamondsuit$}\vskip6pt}

\def\bibitem#1{\par\indent\llap{\rlap{\bf [#1]}\indent}\indent\hangindent
2\parindent\ignorespaces}

\long\def\eatit#1{}

\def\leftheadlinetext{}
\def\rightheadlinetext{}

\def\leftheadline{{\eightrm\folio\hfil \leftheadlinetext\hfil}}
\def\rightheadline{{\eightrm\hfil\rightheadlinetext\hfil\folio}}

\headline={\ifnum\pageno=1\hfil\else
\ifodd\pageno\rightheadline\else\leftheadline\fi\fi}

\def\tenpoint{\def\rm{\fam0\tenrm}
\textfont0=\tenrm \scriptfont0=\sevenrm \scriptscriptfont0=\fiverm
\textfont1=\teni \scriptfont1=\seveni \scriptscriptfont1=\fivei
\def\mit{\fam1} \def\oldstyle{\fam1\teni}
\textfont2=\tensy \scriptfont2=\sevensy \scriptscriptfont2=\fivesy
\def\cal{\fam2}
\textfont3=\tenex \scriptfont3=\tenex \scriptscriptfont3=\tenex
\def\it{\fam\itfam\tenit} 
\textfont\itfam=\tenit
\def\sl{\fam\slfam\tensl} 
\textfont\slfam=\tensl
\def\bf{\fam\bffam\tenbf} 
\textfont\bffam=\tenbf \scriptfont\bffam=\sevenbf
\scriptscriptfont\bffam=\fivebf
\def\tt{\fam\ttfam\tentt} 
\textfont\ttfam=\tentt
\normalbaselineskip=12pt
\setbox\strutbox=\hbox{\vrule height8.5pt depth3.5pt  width0pt}%
\normalbaselines\rm}

\def\eightpoint{\def\rm{\fam0\eightrm}%
\textfont0=\eightrm \scriptfont0=\sixrm \scriptscriptfont0=\fiverm
\textfont1=\eighti \scriptfont1=\sixi \scriptscriptfont1=\fivei
\def\mit{\fam1} \def\oldstyle{\fam1\eighti}%
\textfont2=\eightsy \scriptfont2=\sixsy \scriptscriptfont2=\fivesy
\def\cal{\fam2}%
\textfont3=\tenex \scriptfont3=\tenex \scriptscriptfont3=\tenex
\def\it{\fam\itfam\eightit} 
\textfont\itfam=\eightit
\def\sl{\fam\slfam\eightsl} 
\textfont\slfam=\eightsl
\def\bf{\fam\bffam\eightbf} 
\textfont\bffam=\eightbf \scriptfont\bffam=\sixbf
\scriptscriptfont\bffam=\fivebf
\def\tt{\fam\ttfam\eighttt} 
\textfont\ttfam=\eighttt
\normalbaselineskip=9pt%
\setbox\strutbox=\hbox{\vrule height7pt depth2pt  width0pt}%
\normalbaselines\rm}

\def\largefont{\def\rm{\fam0\largerm}
\textfont0=\largerm \scriptfont0=\largescriptrm \scriptscriptfont0=\tenrm
\textfont1=\largei \scriptfont1=\largescripti \scriptscriptfont1=\teni
\def\mit{\fam1} \def\oldstyle{\fam1\teni}
\textfont2=\largesy 
\def\cal{\fam2}
\def\it{\fam\itfam\largeit} 
\textfont\itfam=\largeit
\def\sl{\fam\slfam\largesl} 
\textfont\slfam=\largesl
\def\bf{\fam\bffam\largebf} 
\textfont\bffam=\largebf 
\scriptscriptfont\bffam=\fivebf
\def\tt{\fam\ttfam\largett} 
\textfont\ttfam=\largett
\normalbaselineskip=17.28pt
\setbox\strutbox=\hbox{\vrule height12.25pt depth5pt  width0pt}%
\normalbaselines\rm}

\def\titlefont{\def\rm{\fam0\titlerm}
\textfont0=\titlerm \scriptfont0=\largescriptrm \scriptscriptfont0=\tenrm
\textfont1=\titlei \scriptfont1=\largescripti \scriptscriptfont1=\teni
\def\mit{\fam1} \def\oldstyle{\fam1\teni}
\textfont2=\titlesy 
\def\cal{\fam2}
\def\it{\fam\itfam\titleit} 
\textfont\itfam=\titleit
\def\sl{\fam\slfam\titlesl} 
\textfont\slfam=\titlesl
\def\bf{\fam\bffam\titlebf} 
\textfont\bffam=\titlebf 
\scriptscriptfont\bffam=\fivebf
\def\tt{\fam\ttfam\titlett} 
\textfont\ttfam=\titlett
\normalbaselineskip=24.8832pt
\setbox\strutbox=\hbox{\vrule height12.25pt depth5pt  width0pt}%
\normalbaselines\rm}

\nopagenumbers

\font\eightrm=cmr8
\font\eighti=cmmi8
\font\eightsy=cmsy8
\font\eightbf=cmbx8
\font\eighttt=cmtt8
\font\eightit=cmti8
\font\eightsl=cmsl8
\font\sixrm=cmr6
\font\sixi=cmmi6
\font\sixsy=cmsy6
\font\sixbf=cmbx6

\font\largerm=cmr12 at 17.28pt
\font\largei=cmmi12 at 17.28pt
\font\largescriptrm=cmr12 at 14.4pt
\font\largescripti=cmmi12 at 14.4pt
\font\largesy=cmsy10 at 17.28pt
\font\largebf=cmbx12 at 17.28pt
\font\largett=cmtt12 at 17.28pt
\font\largeit=cmti12 at 17.28pt
\font\largesl=cmsl12 at 17.28pt

\font\titlerm=cmr12 at 24.8832pt
\font\titlei=cmmi12 at 24.8832pt
\font\titlesy=cmsy10 at 24.8832pt
\font\titlebf=cmbx12 at 24.8832pt
\font\titlett=cmtt12 at 24.8832pt
\font\titleit=cmti12 at 24.8832pt
\font\titlesl=cmsl12 at 24.8832pt

\tenpoint



\def\manyby{\hbox to.75in{\hrulefill}}
\hsize 6.5in 
\vsize 9.2in

\tolerance 3000
\hbadness 3000

\def\item#1{\par\indent\indent\llap{\rlap{#1}\indent}\hangindent
2\parindent\ignorespaces}

\def\itemitem#1{\par\indent\indent
\indent\llap{\rlap{#1}\indent}\hangindent
3\parindent\ignorespaces}

\def\leftheadlinetext{Fitchett, Harbourne \& Holay}
\def\rightheadlinetext{Resolutions of Fat Point Ideals}

\def\mbox#1{\hbox{#1}}
\def\pts#1#2#3{\mbox{${#1}_{#2},\ldots,{#1}_{#3}$}}
\def\pr#1{\mbox{${\bf P}^{#1}$}}
\def\ot{\mbox{$\otimes$}}
\def\h#1#2#3{\mbox{$h^{#1}(#2,#3)$}}
\def\sh#1#2{\mbox{${\cal O}_{#1}(#2)$}}
\def\shn#1{\mbox{${\cal O}_{#1}$}}  
\def\dimcok#1{\mbox{dim(cok($\mu_{#1}$))}} 
\def\dimker#1{\mbox{dim(ker($\mu_{#1}$))}} 

\def\z{\mbox{${\bf Z}$}}
\def\H#1#2#3{\mbox{$H^{#1}(#2,#3)$}}
\def\d#1{\mbox{$#1$}}
\def\cite#1{[\ir{#1}]}
\def\binom#1#2{\hbox{$\left(\matrix{#1\cr #2\cr}\right)$}}

\def\catalisano{C}
\def\curves{F1}
\def\sfmaps{F2}
\def\clsrs{H1}
\def\duke{H2}
\def\mtnwest{H3}
\def\ars{H4}
\def\birmor{H5}
\def\algorithm{H6}
\def\resolution{H7}
\def\igp{H8}
\def\survey{H9}
\def\mumford{M}
\def\nagata{N}
\def\ramanujam{R}
\def\terakawa{T}

\ \vskip-.3in
{\largefont
\hbox to\hsize{\hfil Resolutions of Fat Point Ideals 
involving 8 General Points of \pr2\hfil}
}

\vskip\baselineskip
{\font\authorrm=cmr12 at 14.4pt
\authorrm
\hbox to\hsize{Stephanie Fitchett\hfil }
}
\vskip6pt

\address
Florida Atlantic University, Honors College, Jupiter, FL 33458
email: sfitchet@fau.edu
\endaddress

{\font\authorrm=cmr12 at 14.4pt
\authorrm
\hbox to\hsize{Brian Harbourne${}^*$\hfil }
}
\vskip6pt

\address
Department of Mathematics and Statistics, University of Nebraska-Lincoln
Lincoln, NE 68588-0323
email: bharbour@math.unl.edu
WEB: http://www.math.unl.edu/$\sim$bharbour/
\endaddress

{\font\authorrm=cmr12 at 14.4pt
\authorrm
\hbox to\hsize{Sandeep Holay\hfil }
}
\vskip6pt

\address
Department of Mathematics, Southeast Community College, Lincoln, NE 
email: sholay@sccm.cc.ne.us
\smallskip
November 30, 2000\endaddress
\vskip-\baselineskip

\thanks{\vskip -6pt
\noindent ${}^*$This author benefitted from a National Science 
Foundation grant.
\smallskip
\noindent 1991 {\it Mathematics Subject Classification. } 
MSC: Primary 13P10; 14C99; secondary 13D02; 13H15.
\smallskip
\noindent {\it Key words and phrases. }  Minimal free resolution, 
rational surface, fat point.\smallskip}

\vskip\baselineskip
\Abstract{Abstract: 
The main result, \ir{mainthmgeom}, 
provides an algorithm for determining the
minimal free resolution of fat point subschemes of \pr2 involving
up to 8 general points of arbitrary multiplicities; the resolutions
obtained hold for any algebraically closed field,
independent of the characteristic. 
The algorithm works by giving a formula in nice cases, 
and a reduction to the nice cases otherwise. 
The algorithm, which does not involve 
Gr\"obner bases, is very fast. Partial information is also 
obtained in certain cases
with $n>8$.}
\vskip\baselineskip

\irrnSection{Introduction}{intro}
Determining the Hilbert function and minimal free resolution
of ideals defining $n$ general fat points of \pr2 is a difficult
problem that has attracted the attention of numerous researchers
over the years. For $n>9$, the problem remains unsolved in general.
For $n\le9$, Nagata \cite{nagata} resolved the problem 
for Hilbert functions, while for $n\le5$, Catalisano \cite{catalisano}
resolved the problem for minimal free resolutions, extended
(using different methods) by Fitchett \cite{sfmaps} to
$n=6$ and by Harbourne \cite{algorithm} to $n=7$. In this paper,
we now extend this work to $n=8$. Our results are based on studying
linear systems on blow ups of \pr2 at the $n$ points.
For $n\le 8$, these surfaces are Del Pezzo surfaces,
and hence the semigroup of classes of effective divisors
is finitely generated. It will be very difficult to extend 
our approach to the case
of $n=9$; for one thing, the surfaces obtained for $n\ge9$
are no longer Del Pezzos, for another, the semigroup
of effective divisor classes is no longer finitely generated.

We now recall the notion of fat points.
Consider  $n$ distinct points \pts{p}{1}{n} of \pr2. Given nonnegative 
integers $m_i$, a fat point subscheme $Z=m_1p_1+\cdots+m_np_n$ is that 
subscheme defined by the homogeneous ideal 
$I_Z=I_1^{m_1}\cap\cdots\cap I_n^{m_n}$
in the homogeneous coordinate ring $R=k[\pr2]$ of \pr2
(over any algebraically closed field $k$), where $I_j$
is the ideal generated by all forms vanishing at $p_j$.

Because $Z$ has codimension 2 and is arithmetically Cohen-Macaulay,
the minimal free graded resolution of $I_Z$ is of the form
$0\to F_1\to F_0\to I_Z\to 0$. To determine $F_0$ it then suffices to
compute for each $t$ the dimension $\nu_{t+1}(Z)$ of the 
cokernel of the multiplication map
$\mu_t(Z):(I_Z)_t\ot R_1\to (I_Z)_{t+1}$, since $F_0$ is
simply $\oplus_{i\geq 0} R[-i]^{\nu_i}$. 
If we know the 
Hilbert function $h_Z$ of $I_Z$
(i.e., the dimension $h_Z(t)=\hbox{dim}(I_Z)_t$ 
for each degree $t$ of the homogeneous component 
$(I_Z)_t$ of $I_Z$ of degree $t$), exactness of the resolution
then allows us to determine the Hilbert function of $F_1$ and thus 
(since $F_1$ is free) $F_1$ itself. More explicitly, since
in general $\nu_t(Z)=\hbox{dim}(\hbox{Tor}_0(I_Z,k))_t$ and
$F_1=\oplus_{i\geq 0} R[-i]^{s_i}$ where 
$s_i=\hbox{dim}(\hbox{Tor}_1(I_Z,k))_i$, we can tensor
the Koszul complex for $k$ by $I_Z$ to compute
these Tor's, and we find $\nu_i-s_i=\Delta^3h_Z(i)$,
where $\Delta$ is the difference operator (hence 
$\Delta h_Z(i)=h_Z(i)-h_Z(i-1)$, for example). 

Instead of working with $I$ and its components, 
by a standard translation, one can instead work on the surface $X$ 
obtained by blowing up the $n$ distinct points \pts{p}{1}{n} of \pr2, 
in which case one then has the corresponding birational 
morphism $\pi:X\to\pr2$. Throughout this paper, this is 
what $X$ will denote.

We now recall the standard translation referred to above.
With respect to $\pi: X\to\pr2$, let $L$ be the 
total transform to $X$ of a line 
on \pr2, and let $E_i = \pi^{-1}(p_i)$ for $1 \leq i \leq n$.  
We call the divisors $L,  E_1, \ldots ,  E_n$ an {\it exceptional
configuration};
their classes $[L],  [E_1], \ldots ,  [E_n]\in \hbox{Cl}(X)$ form a basis
for the divisor class group $\hbox{Cl}(X)$ of $X$, with intersections given 
by $E_i \cdot  E_j = - \delta_{ij}$ (where $\delta_{ij}$ denotes
Kronecker's delta), $L \cdot E_i = 0$ for all $i$, 
and $L^2 = 1$. For $Z=m_1p_1 + \cdots + m_np_n$, under natural identifications
we have $(I_Z)_t =\H{0}{X}{\sh{X}{F_t}}$, 
hence $h_Z(t)=\h{0}{X}{\sh{X}{F_t}}$, where 
$F_t =  tL-m_1E_1 - \cdots -  m_nE_n$.  Moreover, the map 
$\mu_t(Z): (I_Z)_t \ot R_1 \to (I_Z)_{t+1}$ given
by multiplication corresponds under these identifications to
the natural map
$\mu_{F_t}:  \H{0}{X}{\sh{X}{F_t}} \ot \H{0}{X}{\sh{X}{L}} 
\to \H{0}{X}{\sh{X}{\sh{X}{F_t}\ot \sh{X}{L}}} = \H{0}{X}{\sh{X}{F_{t+1}}}$.
Thus, although in this paper we will be concerned with $\mu_F$
for an arbitrary divisor $F$,
our results have an immediate application to computing resolutions of 
fat point ideals. We also note that 
as long as $\h{0}{X}{\sh{X}{H}}$ is known for arbitrary divisors $H$,
to compute the dimension of the cokernel of $\mu_F$ for some divisor
$F$, it is just as good to compute the dimension of the kernel of $\mu_F$
or the rank of $\mu_F$, depending on convenience.

In addition to assuming that we have a birational
morphism $X\to\pr2$ (i.e., that $X$ is {\it basic}, that is, obtained
by blowing up $n\ge0$ points $p_i$ of \pr2), we will assume
that there is a particular smooth, irreducible anticanonical
divisor $D_X$ on $X$ (i.e., the points blown up lie on a smooth
cubic), and that $X$ satisfies the following
properties:

\item{$\bullet$} (A1) the only integral curves besides possibly $D_X$ of negative 
self-intersection on $X$ are exceptional curves (i.e., 
smooth rational curves of self-intersection $-1$), and 
\item{$\bullet$} (A2) $\h{1}{X}{\sh{X}{F}}=0$ for any effective,
numerically effective divisor $F$ 
($F$ being {\it numerically effective} means 
that $F\cdot H\geq 0$ for every effective divisor $H$).

\noindent In these circumstances we will say that $X$ 
is a {\it good} surface.

For example, any blow up $X$ of \pr2 at $n\le 8$ 
general points is good: (A2) holds by Theorem 8 of
\cite{mtnwest}, and (A1) holds by adjunction since
the anticanonical class 
$-K_X$ is ample (see the first two paragraphs of
the proof of Theorem 1 of \cite{mtnwest}).
For another example, say that a basic
surface $X$ with a
smooth, irreducible anticanonical
divisor $D_X$ is {\it injective} if the canonical map
$\hbox{Pic}(X)\to\hbox{Pic}(D_X)$ is injective.
(Although we will sometimes use $D$ to denote an
arbitrary divisor, we will throughout this paper use $D_X$ to 
denote a smooth, irreducible anticanonical divisor.)
Such a surface $X$ is good:
(A2) holds by Theorem III.1 of \cite{ars},
while (A1) holds by adjunction (any integral curve $C$ of
negative self-intersection which is neither an exceptional curve
nor $D_X$ must by adjunction have $C\cdot D_X=0$ and hence $C=0$ since
$\hbox{Pic}(X)\to\hbox{Pic}(D_X)$ is injective).

We note that
the condition that $\hbox{Pic}(X)\to\hbox{Pic}(D_X)$ is injective
holds in all characteristics if the points blown up
are sufficiently general points of a smooth plane cubic curve, 
if the ground field $k$
is sufficiently large. It cannot hold, however, if 
$k$ is the algebraic closure of a finite field.
On the other hand, assuming the points $p_i$ are general
points of a smooth plane cubic curve,
our result for the rank of
$\mu_F$ for a particular $F$ or for the minimal free resolution
of the ideal $I_Z$ for a particular $Z$, holds over any
algebraically closed field $k$ (this is because our result 
for a particular $F$ or $Z$ holds 
for some set of $n$ $\bar K$-points of a smooth plane cubic $k$-curve
$C\subset\pr2$, where $\bar K$ is the algebraic closure of
a sufficiently large extension $K$ of $k$, but the
$k$-points of $C^n$ are a dense subset of the $\bar K$-points).

We will denote the set of reduced, irreducible curves $C$ on $X$ 
with $\h{0}{X}{\sh{X}{C}}=1$ by $\Gamma_X$. By Theorem III.1 of \cite{ars},
(A1) implies that $\Gamma_X$ includes the set of exceptional curves
on $X$, and $D_X$ if $D_X^2\le0$, but nothing else.

Moreover, as shown 
in section I.4 of \cite{algorithm}
({\it mutatis mutandis}, since $-K_X$ need not 
always be the class of an effective divisor 
in \cite{algorithm}), given that $X$ is good,
one can in a completely effective manner 
determine the fixed components and dimension 
of any complete linear system on $X$. 
Our main interest, then, is, in addition, to
determine the rank of $\mu_F$ for an arbitrary divisor $F$.

Assuming (A1) and (A2) above, our approach is
to reduce the problem of computing the rank of $\mu_F$ 
for an arbitrary divisor $F$ to that of doing so for a special class
of divisors. A similar approach was taken in \cite{algorithm}, in which
a reduction was made to ample $F$. Unfortunately,
the class of ample divisors is still too coarse to 
get a nice answer in general (note that the exceptional case
in \ir{mainthmgeom}(c)(ii) is ample), even for $n=8$ points,
so here we present a refined reduction to a particular special
set of divisors $F$. Although for arbitrary $n$ we cannot always handle the 
resulting divisors by present methods, for $n\leq 8$ 
using ad hoc methods we always can, thereby achieving
a complete determination. In terms of fat points, this
gives a complete determination of the minimal free resolution
for the ideal of any fat points subscheme involving up to 8 general
points of \pr2, essentially by reduction to nice cases
where a formula applies.

To state our main result, we introduce some terminology and notation.
For convenience, we will often write simply
$\h{0}{X}{F}$ for $\h{0}{X}{\sh{X}{F}}$.
Also, we say a divisor $F$ on $X$ is {\it monotone} provided 
$F \cdot E_1 \geq F \cdot E_2 \geq \cdots \geq F \cdot E_n$.
We now define quantities $\lambda$ and $\Lambda$ for each curve $C$ 
in $\Gamma_X$.
For $C=E_i$ for any $i$, let $\lambda_C=\Lambda_C=0$.
Otherwise, let $m_C$ be the maximum of $C\cdot E_1, \ldots, C\cdot E_n$,
define $\Lambda_C$ to be the maximum
of $m_C$ and of $(C \cdot L) - m_C$, define
$\lambda'_C$ to be the minimum of $m_C$ and of $(C \cdot L) - m_C$, 
and, if $C$ is a smooth rational curve, define
$\lambda_C$ to be $\lambda'_C$, and otherwise
define $\lambda_C$ to be the maximum of $\lambda'_C$ and 2.
For example, if $C$ is anticanonical,
then $\lambda'_C=1$ and $\lambda_C=\Lambda_C=2$.

Here now is our main result, \ir{mainthmgeom}. 
For general $n$, it determines the rank of $\mu_F$ for
certain (sufficiently nice) divisors $F$, which is sometimes
enough to determine resolutions, as shown in \ir{sectionexample}.
For $n\leq8$, \ir{mainthmgeom}
gives a reduction for arbitrary $F$ to the 
nice case, thereby in principle providing an algorithm
for computing the rank of $\mu_F$ for any divisor $F$. 
Briefly, the algorithm works as follows. By applying \ir{eff}
with \ir{lists}(a), we can compute $\h{0}{X}{F}$
for any divisor $F$ (see section 4 for 
examples, or \cite{clsrs}, \cite{ars} and \cite{algorithm} for more 
comprehensive background). 
By \ir{lists}, $\Gamma_X$ is finite for $n=8$, 
so cases (a) and (b) of \ir{mainthmgeom} 
are easy to implement. If cases (a) or (c) occur, we are done,
while case (b) need be applied only if $\h{0}{X}{F}>0$
since $\mu_F$ is clearly injective if $\h{0}{X}{F}=0$.
In case (b) with $\h{0}{X}{F}>0$,
repeated applications of (b) gives a divisor $F$
for which either $\h{0}{X}{F}=0$
or which satisfies the conditions either of (a) or (c).
(It is not hard to implement the 
algorithm explicitly; see \cite{survey}, which includes 
an explicit Macaulay script implementing this algorithm.)

\prclm{Theorem}{mainthmgeom}{Let $X$ be a good surface 
with exceptional configuration
$L, E_1,\ldots,E_n$. Let $F$ be a monotone divisor on $X$.  
\item{(a)} If $F \cdot C \geq \Lambda_C$ 
for all $C\in \Gamma_X$, then $\mu_F$ has {\rm maximal rank}
(i.e., is either injective or surjective), hence
$\dimcok{F}=\hbox{max}(0, \h{0}{X}{F+L}-3\h{0}{X}{F})$.
\item{(b)} If $F \cdot C < \lambda_C$ for some $C\in \Gamma_X$, 
then $\dimker{F}=\dimker{F-C}.$
\item{(c)} If $n=8$ but neither case (a) nor case (b) holds, 
then either
\itemitem{(i)} $F \cdot (L - E_1 - E_2)=0$, in which case 
the dimension of $\mbox{cok}(\mu_F)$ is
$\h{1}{X}{F - (L - E_1)} + \h{1}{X}{F - (L - E_2)}$, or
\itemitem{(ii)} $[F]$ is $[3L - E_1 - \cdots - E_7] + 
r[8 L - 3E_1 - \cdots - 3E_7 - E_8]$ for some 
$r \geq 1$, in which case  $\dimcok{F} = r$ and 
$\dimker{F} = r+1$, or
\itemitem{(iii)} $\mu_F$ has maximal rank, hence 
$\dimcok{F}=\hbox{max}(0,F\cdot L+F\cdot K_X-F^2)$.}

\Prf 
Part (a) is \ir{nicecasecor}. (The proof 
here depends on a criterion for $\mu_F$ to be surjective 
which involves vanishing of certain $h^1$'s, in the form of 
$q^*(F)$ and $l^*(F)$ defined below; see \ir{qstarlstarlemma}(b). 
Being defined in terms of $h^1$'s involving $F$, 
intuitively we can expect $q^*(F)$ and $l^*(F)$ to vanish if
$F$ is sufficiently nice, which is precisely what our hypotheses
guarantee.)   
Part (b) is \ir{samekernel}. (This part generalizes the
fact that $\mu_F$ and $\mu_{F-C}$ have kernels of the same 
dimension if $F$ is effective and $C$ is an integral curve 
with $F\cdot C<0$, which is obvious since in this situation
$C$ is a fixed component of $|F|$.)
Most of the work here is in proving (c),
which follows from \ir{partcprop} and 
\ir{specialfamily}. (The proof of this part
amounts to an analysis of all possibilities not dealt with
by (a) or (b).)
\qed

\irrnSection{The Nice Case and the Reduction}{nicecase}
We regard case (a) of \ir{mainthmgeom} as the nice case because
there we have the rank of $\mu_F$ directly. Case (b) of
\ir{mainthmgeom} is the reduction case.
In this section we prove cases (a) and (b) of \ir{mainthmgeom}.
Given a basic surface $X$ with exceptional configuration
$L, E_1,\ldots, E_n$ and a monotone divisor $F$ on $X$, we will
denote \h{0}{X}{F-E_1} by $q(F)$, \h{1}{X}{F-E_1} by $q^*(F)$,
\h{0}{X}{F-(L-E_1)} by $l(F)$ and \h{1}{X}{F-(L-E_1)} by $l^*(F)$.

\prclm{Lemma}{qstarlstarlemma}{Let $X$ 
be obtained by blowing up distinct points of \pr2,
with exceptional configuration $L, E_1,\ldots,E_n$.
Let $F$ be a monotone divisor on $X$.
\item{(a)} Then $\dimker{F}\le q(F)+l(F)$; in particular,
$\mu_F$ is injective if $q(F)=0=l(F)$.
\item{(b)} If $F$ is effective and $\h{1}{X}{F}=0$, 
then $\dimcok{F}\le q^*(F)+l^*(F)$; in particular, 
$\mu_F$ is surjective, if, in addition, $l^*(F)=0=q^*(F)$.}

\Prf
(a) See Lemma 4.1 of \cite{igp} (which assumes
$F\cdot E_i>0$ for all $i$, but it is easy to check
that the result holds even if $F\cdot E_i\le 0$ for some $i$).

(b) Clearly, $L$ is numerically effective.
Thus, $F\cdot L\ge 0$, since $F$ is effective.
Now, from $\h{1}{X}{F}=0$
and $0\to \sh{X}{F}\to \sh{X}{F+L}\to \sh{L}{F+L}\to 0$,
we see $\h{1}{X}{F+L}$ vanishes also and we compute
$\h{0}{X}{F+L}-3\h{0}{X}{F}=2+F\cdot L-2\h{0}{X}{F}$.
Similarly, $l^*(F)-l(F)=F\cdot (L-E_1) + 1-\h{0}{X}{F}$
and $q^*(F)-q(F)=F\cdot E_1 + 1-\h{0}{X}{F}$,
so $(l^*(F)-l(F))+(q^*(F)-q(F))=\h{0}{X}{F+L}-3\h{0}{X}{F}$.
Therefore,
$\dimcok{F}=\dimker{F}+
\h{0}{X}{F+L}-3\h{0}{X}{F}\le l(F)+q(F)+\h{0}{X}{F+L}-3\h{0}{X}{F}
= l^*(F)+q^*(F)$, as required. 
\qed

We will need to refer to the following result.

\prclm{Lemma}{lists}{Let $C$ be a curve on the blow-up $X$ of \pr2 at 8
general points, with exceptional configuration
$L, E_1,\ldots,E_8$.
\item{(a)} Then, up to permutation of the indices, 
$C$ is an exceptional curve if and only if 
$[C]$ is one of the following: 
$[E_8]$, $[L-E_1-E_2]$, $[2L-E_1-\cdots-E_5]$,
$[3L-2E_1-E_2-\cdots-E_7]$, $[4L-2E_1-2E_2-2E_3-E_4-\cdots-E_8]$,
$[5L-2E_1-\cdots-2E_6-E_7-E_8]$, or $[6L-3E_1-2E_2-\cdots-2E_8]$.
\item{(b)} And, up to permutation of the indices, 
$C$ is a smooth rational curve with $C^2 = 0$ if and only if 
$[C]$ is one of the following: 
$[L-E_1]$, 
$[2L-E_1-\cdots-E_4]$,
$[3L-2E_1-E_2-\cdots-E_6]$, 
$[4L-2E_1-2E_2-2E_3-E_4-\cdots-E_7]$,
$[4L-3E_1-E_2-\cdots-E_8]$, 
$[5L-3E_1-2E_2-2E_3-2E_4-E_5-\cdots-E_8]$, 
$[5L-2E_1-\cdots-2E_6-E_7]$, 
$[6L-3E_1-3E_2-2E_3-\cdots-2E_6-E_7-E_8]$,
$[7L-3E_1-\cdots-3E_4-2E_5-2E_6-2E_7-E_8]$, 
$[7L-4E_1-3E_2-2E_3-\cdots-2E_8]$,
$[8L-3E_1-\cdots-3E_7-E_8]$, 
$[8L-4E_1-3E_2-\cdots-3E_5-2E_6-2E_7-2E_8]$,
$[9L-4E_1-4E_2-3E_3-\cdots-3E_7-2E_8]$, 
$[10L-4E_1-\cdots-4E_4-3E_5-\cdots-3E_8]$, and 
$[11L-4E_1-\cdots-4E_7-3E_8]$.}

\Prf
Under the action on $\hbox{Cl}(X)$ by the Weyl group
(which is generated by permutation of the indices 
and by the action of quadratic
Cremona transformations centered at any three of the 8 points;
see \cite{clsrs} or \cite{duke}),
any class of an effective divisor $F$ is in the orbit of
a class $F'$ of the form $[dL+a_1E_1+\cdots+a_8E_8]$, where
$d\ge 0$, $3d+a_1+\cdots+a_8\ge 0$ (since, as mentioned
in the introduction, $-K_X$ is ample),
$d+a_1+a_2+a_3\ge 0$ and $a_1\le a_2\le \cdots\le a_8$.
If $F$ is reduced and irreducible, then 
$F'$ is either $[E_8]$ or $a_8\le 0$. In the latter case,
$F'$ is by Lemma 1.4 of \cite{clsrs} a nonnegative sum of the classes
of $L$, $L-E_1$, $2L-E_1-E_2$, $3L-E_1-E_2-E_3$, $\ldots$, 
$3L-E_1-\cdots-E_8$. But any such sum
has nonnegative self-intersection. Thus classes of exceptional curves
are in the orbit of $[E_8]$; these are the classes listed in (a).
And the only such sums with self-intersection 0 are the positive
multiples of $[L-E_1]$, hence only $[L-E_1]$ itself represents the class of
a reduced and irreducible divisor with self-intersection 0. 
Thus the classes to be listed in (b)
comprise the orbit of $[L-E_1]$ under the action of the Weyl group;
it is easy to check that the list in (b)
is (up to permutations) the complete orbit.
\qed

Given a surface $X$, denote by EFF${}_X$ (or just EFF) the 
subsemigroup of the divisor class group
$\mbox{Cl}(X)$ of $X$ of classes
of effective divisors, and let NEFF denote the cone of 
classes of numerically effective divisors.

\prclm{Lemma}{eff}{Let $X$ be a good surface
obtained by blowing up $n$ points of \pr2. 
\item{(a)} If $2\le n<8$, then EFF is generated by the classes of
exceptional curves, while for $n\geq8$, EFF is generated
by the classes of exceptional curves and by $-K_X$.
\item{(b)} Let $F\in\hbox{Cl}(X)$. 
If $2\le n\leq8$, then $F$ is in NEFF if and only
if $F\cdot C\geq0$ for every exceptional curve $C$, while
if $9\le n$, then $F$ is in NEFF if and only
if both $-F\cdot K_X\geq0$ and $F\cdot C\geq0$ for 
every exceptional curve $C$.
\item{(c)} EFF contains NEFF, and
$\h{1}{X}{F}=\h{2}{X}{F}=0$  
and $\h{0}{X}{F}=(F^2-F\cdot K_X)/2+1$ for any $F\in\hbox{NEFF}$.}

\Prf
(a) By (A1), any effective divisor is (up to linear equivalence) a nonnegative
sum of $D_X$, exceptional curves and an element of NEFF.
But by Corollary 3.2 and Lemma 1.4, both of \cite{clsrs},
any element of NEFF is (with respect to some exceptional configuration
$L$, $E_1,\ldots,E_n$)
a sum of $[L]$, $[L-E_1]$, $[2L-E_1-E_2]$, $-K_X$ and classes of exceptional
curves. But $n\ge 2$, so $[L]=[L-E_1-E_2]+[E_1]+[E_2]$,
$[L-E_1]=[L-E_1-E_2]+[E_2]$ and $[2L-E_1-E_2]=2[L-E_1-E_2]+[E_1]+[E_2]$
are sums of classes of exceptional curves, so NEFF
is generated by the classes
of exceptional curves and by the class $-K_X$ of $D_X$.
Moreover, for $n=2$, 3, 4, 5, 6 or 7,
$-K_X$ is, respectively, the following sums of classes
of exceptional curves: $3[L-E_1-E_2]+2[E_1]+2[E_2]$,
$[L-E_1-E_2]+[L-E_1-E_3]+[L-E_2-E_3]+[E_1]+[E_2]+[E_3]$,
$[L-E_1-E_2]+[L-E_3-E_4]+[L-E_1-E_2]+[E_1]+[E_2]$, 
$[L-E_1-E_2]+[L-E_3-E_4]+[L-E_1-E_5]+[E_1]$,
$[L-E_1-E_2]+[L-E_3-E_4]+[L-E_5-E_6]$, or
$[2L-E_1-E_2-\cdots-E_5]+[L-E_6-E_7]$.

(b) This follows immediately from (a),
except in the case that $n=8$. For $n=8$, we have
$-2K_X=[3L-2E_1-E_2-\cdots-E_7]+[3L-E_2-\cdots-E_7-2E_8]$,
hence if $F\cdot E\geq0$ for all exceptional curves $E$, then also
$-F\cdot K_X\geq0$, and our result follows here too.

(c) By Proposition 4 of \cite{mtnwest}, $F^2\ge0$ and $\h{2}{X}{F}=0$ for
any $F\in \hbox{NEFF}$. Since we assume $-K_X$ is effective,
we also have $-F\cdot K_X\ge 0$. But by Riemann-Roch,
$\h{0}{X}{F}\ge (F^2-F\cdot K_X)/2+1$, so $\h{0}{X}{F}\ge 1$,
hence $F\in \hbox{EFF}$. Now
$\h{1}{X}{F}=0$ for any $F\in\hbox{NEFF}$ by assumption (A2)
and the rest is immediate from Riemann-Roch.
\qed

Here is the proof of \ir{mainthmgeom}(a).

\prclm{Corollary}{nicecasecor}{Let $X$ be a good 
surface with exceptional configuration
$L,  E_1, \ldots,  E_n$, and let $F$ be a monotone divisor on $X$. 
If $F \cdot C \geq \Lambda_C$ 
for all $C\in \Gamma_X$, then $\mu_F$ has maximal rank.}

\Prf
If $\h{0}{X}{F} = 0$, then clearly $\mu_F$ is injective
and so has maximal rank.
So assume $\h{0}{X}{F}>0$.
Since both $\Lambda_C\ge0$ and
$F\cdot C\ge\Lambda_C$ for all $C\in \Gamma_X$  and hence
for all curves $C$ of negative self-intersection,
we see $F$ is numerically effective.

If $n\le 5$, then $\mu_F$ is surjective simply because
$F$ is numerically effective (Theorem 3.1.2, \cite{resolution}).
So we may assume that $n\ge 6$. 

Since $F$ is numerically effective, then $\h{1}{X}{F}=0$
by (A2) and $F\cdot E_i\ge 0$ for all $i$, and so
$(F-E_1)\cdot E_i\ge 0$ for all $i$. 
If $C$ is an exceptional curve but not $E_i$ for any $i$,
then $\Lambda_C\geq C\cdot E_1$ and hence $F\cdot C\geq \Lambda_C$
implies $(F-E_1)\cdot C\geq 0$. If $6\le n\le 8$, then
$F-E_1$ is numerically effective by \ir{eff}, hence
$q^*(F)=\h{1}{X}{F-E_1}=0$ by (A2). 
If $n\geq9$, then $D_X\in\Gamma_X$, so in addition we have
$(F-E_1)\cdot D_X\geq\Lambda_{D_X}-1\geq1$, and again we see
$F-E_1$ is numerically effective, and 
$q^*(F)=\h{1}{X}{F-E_1}=0$. 

Now consider $l^*(F)$. If $(F-(L-E_1))\cdot E_i<0$ for some
$i$, then clearly $i=1$ and $F\cdot E_1=0$. Since $F$ is monotone
we see $F\cdot E_i=0$ for all $i$, so $F$ is (up to
linear equivalence) a nonnegative multiple of $L$, 
for which it is easy to
see that $\mu_F$ has maximal rank. Thus we may assume 
$(F-(L-E_1))\cdot E_i\ge 0$ for all $i$. For any other 
exceptional curve $C$ we have $F\cdot C\ge\Lambda_C$,
so $(F-(L-E_1))\cdot C\ge\Lambda_C-(L-E_1)\cdot C\ge0$
too. Thus, as above, $F-(L-E_1)$ is numerically effective if
$n\leq8$, so $l^*(F)=\h{1}{X}{F-(L-E_1)}=0$, while 
if $n\ge 9$, we have in addition that
$(F-(L-E_1))\cdot D_X\geq\Lambda_{D_X}-2\geq0$
so again $l^*(F)=\h{1}{X}{F-(L-E_1)}=0$.
The result now follows by \ir{qstarlstarlemma}(b).
\qed

We now prove part (b) of \ir{mainthmgeom}.

\prclm{Proposition}{samekernel}{Let $F$ be a divisor on a good surface $X$
having exceptional configuration $L$, $E_1,\ldots,E_n$.
If $F \cdot C < \lambda_C$ for some $C\in \Gamma_X$, 
then $\dimker{F}=\dimker{F-C}$.}

\Prf
If $\h{0}{X}{\sh{X}{F}}=0$, then, of course, 
neither \sh{X}{F} nor \sh{X}{F-C}
have any global sections, so
both $\mbox{ker}(\mu_F)$ and $\mbox{ker}(\mu_{F-C})$
vanish. Therefore, we may assume $F$ is effective.

If $C$ is a fixed component of $|F|$ 
for some $C\in\Gamma_X$, then clearly
the canonical injection $\sh{X}{F-C}\to\sh{X}{F}$
induces an isomorphism both of global sections
and of kernels of $\mu$.
So now we may assume that $|F|$ is fixed component free,
hence $F\cdot C\geq0$ for all $C\in\Gamma_X$.

Since $\lambda_{D_X}=2$, if $F\cdot D_X<\lambda_{D_X}$, then 
it must be that $F\cdot D_X$ is 0 or 1. If 0, then \sh{X}{F} is in the kernel
of $\hbox{Pic}(X)\to\hbox{Pic}(D_X)$, and since this is injective,
we see $F=0$, in which case both
$\mbox{ker}(\mu_F)$ and $\mbox{ker}(\mu_{F-D_X})$ again vanish.
If 1, then $\h{0}{D_X}{\sh{D_X}{F}}=1$ so
$\mu_{\sh{D_X}{F}}:\H{0}{D_X}{\sh{D_X}{F}} \ot \H{0}{X}{L} 
\to \H{0}{D_X}{\sh{D_X}{F+L}}$ has the same kernel as the map
$\H{0}{X}{L} \to \H{0}{D_X}{\sh{D_X}{L}}$ given by restriction,
which is easily seen to be injective.
From the exact sequence
$0\to\mbox{ker}(\mu_{F-D_X})\to\mbox{ker}(\mu_{F})\to
\mbox{ker}(\mu_{\sh{D_X}{F}})$
induced by $0\to \sh{X}{F-D_X}\to\sh{X}{F}\to \sh{D_X}{F}\to 0$ 
we see that $\dimker{F} = \dimker{F-D_X}$.

Finally, if $F \cdot C < \lambda_C$ for some $C\in \Gamma_X$
and $C\neq D_X$, then $C$ is an exceptional curve.
As above, the exact sequence $0\to \sh{X}{F-C}\to\sh{X}{F}\to \sh{C}{F}\to 0$ 
induces an exact sequence
$0\to\mbox{ker}(\mu_{F-C})\to\mbox{ker}(\mu_{F})\to \mbox{ker}(\mu_{F\cdot
C,C})$, 
where $\mu_{F\cdot C,C}$ denotes the map 
$$\H{0}{C}{\sh{X}{F}\ot\shn{C}}\ot\H{0}{X}{L}\to\H{0}{C}{\sh{X}{F+L}\ot\shn{
C}}.$$
By Theorem 3.1 of \cite{curves}, $\mbox{ker} \mu_{F\cdot C,C}$ 
is injective if $F \cdot C < \lambda_C$, so $\dimker{F}$ equals
$\dimker{F-C}$.
\qed

\irrnSection{The Main Theorem, Part (c)}{partc}
To prove Part (c) of \ir{mainthmgeom}, we need some additional
background, which we now develop. For the purpose
of stating the next result, given sheaves $\cal F$ and $\cal L$
on a scheme $Y$, 
we denote the kernel of  
$$\H{0}{Y}{{\cal F}} \ot \H{0}{Y}{{\cal L}} 
\to \H{0}{Y}{{\cal F}\ot {\cal L}}$$
by $R({\cal F}, {\cal L})$ and the cokernel by $S({\cal F}, {\cal L})$
(taking $Y$ to be understood).

\prclm{Proposition}{mumes}{Let $Y$ be a closed subscheme of projective space, let 
${\cal F}$ and ${\cal L}$
be coherent sheaves on $Y$ and let ${\cal C}$ be the 
sheaf associated to an effective Cartier divisor $C$
on $Y$.  If the restriction homomorphisms $\H{0}{Y}{\cal F}
\to \H{0}{C}{{\cal F} \ot {\cal O}_C}$ and $\H{0}{Y}{{\cal F}
\ot {\cal L}} \to \H{0}{C}{{\cal F} \ot {\cal L} \ot {\cal O}_C}$
are surjective, then we have an exact sequence
$$\eqalign{0\to & R({\cal F} \ot {\cal C}^{-1}, {\cal L})\to R({\cal F}, {\cal L})\to 
R({\cal F} \ot {\cal O}_C, {\cal L})\to\cr
                & S({\cal F} \ot {\cal C}^{-1}, {\cal L})\to S({\cal F}, {\cal L})\to 
S({\cal F} \ot {\cal O}_C, {\cal L})\to0.\cr}$$}

\Prf
This is a snake lemma argument; see \cite{mumford}.
\qed

\rem{Example}{anexmpl}{The preceding result will 
often be applied in situations
in which we have an exact
sequence $0\to \sh{X}{F-C}\to\sh{X}{F}\to\sh{C}{F}\to 0$,
where $C$ is a curve on a surface $X$ with an exceptional configuration
$L, E_1,\ldots,E_n$, and $F$ is a divisor on $X$ with
$\h{1}{X}{F-C}=0$ and $\h{1}{L}{F-C+L}=0$. 
In such a situation, we can apply
\ir{mumes} with $Y=X$, ${\cal F}=\sh{X}{F}$,
${\cal C}=\sh{X}{C}$ and ${\cal L}=\sh{X}{L}$. For example,
let $X$ be a blow up of \pr2 at 8 general points, with exceptional 
configuration $L, E_1, \ldots, E_8$. Recall $D_X$ is a
smooth, irreducible anticanonical divisor on $X$. For future 
reference we would like to show that $\mu_{D_X+tL}$ and $\mu_{2D_X+tL}$
have maximal rank for all integers $t$. First say that $t=0$.
Consider the exact sequence 
$0\to \sh{X}{2D_X-E}\to\sh{X}{2D_X}\to\sh{E}{2}\to 0$,
where $E$ is the exceptional curve whose class is
$[6L-3E_1-2E_2-\cdots-2E_8]$. Since $2D_X-E=E_1$,
clearly $\mu_{2D_X-E}$ is injective, and
by \ir{squareZero} below, $\mu_{2,E}$
is also injective, hence applying \ir{mumes}
it follows $\dimker{2D_X}=0$; i.e., $\mu_{2D_X}$
has maximal rank, and an easy calculation
now shows that $\dimcok{2D_X}=0$ too.
From \ir{mumes} and the exact sequence 
$0\to \sh{X}{D_X}\to\sh{X}{2D_X}\to\sh{D_X}{2D_X}\to 0$
it also follows that $\dimker{D_X}=0$.
It is easier to see that $\mu_{D_X+tL}$ and
$\mu_{2D_X+tL}$ have maximal rank for all $t\ne 0$,
since $\h{0}{X}{D_X+tL}=0$ and 
$\h{0}{X}{2D_X+tL}=0$ for $t<0$
(in which case $\mu_{D_X+tL}$ and
$\mu_{2D_X+tL}$ are clearly injective),
while for $t\ge0$ we have $\h{1}{X}{D_X+tL}=0$
and $\h{1}{X}{2D_X+tL}=0$, in which case we apply
the general fact that $\mu_{H+L}$ is surjective
if $\h{1}{X}{H+tL}=0$ for $t\ge0$. (Indeed, $\h{1}{X}{H}=0$ ensures
$\H{0}{X}{H+L}\to\H{0}{L}{H+L}$ is surjective,
and hence $S(\sh{L}{L},\sh{X}{H+L})=S(\sh{L}{L},\sh{L}{H+L})$. Then
applying \ir{mumes} with $Y=X$, $C=L$,
${\cal F}=\sh{X}{L}$ and ${\cal L}=\sh{X}{H+L}$
we see $\mu_{H+L}$ is onto---just note that
$S(\sh{X}{L}, \sh{X}{H+L})=0$ since in this situation it is easy
to check that both 
$S(\shn{X}, \sh{X}{H+L})=0$ and
$S(\sh{L}{L},\sh{L}{H+L})=0$.)}

\prclm{Proposition}{perptodegone}{Let $F$ be a monotone divisor on 
a good surface $X$ with exceptional configuration
$L,  E_1, \ldots ,  E_n$.
If  $F \cdot (L - E_1 - E_2) = 0$, then 
$$\dimker{F} = \h{0}{X}{F - (L - E_1)} + \h{0}{X}{F - (L - E_2)}.$$
If in addition $F$ is effective and $\h{1}{X}{F}=0$, then
$$\dimcok{F} = \h{1}{X}{F - (L - E_1)} + \h{0}{X}{F - (L - E_2)}.$$}

\Prf
If $\h{0}{X}{F}=0$, then certainly $\dimker{F}=0$ and both
$\h{0}{X}{F - (L - E_1)}$ and $\h{0}{X}{F - (L - E_2)}$ are also 0. Otherwise,
this is Proposition II.2(e) and Remark II.3, both of \cite{algorithm}.
\qed

We recall that a nonzero element of a free abelian group is {\it primitive}
if it is not a multiple greater than 1 of another element of the group.

\prclm{Lemma}{hOneisZero}{Let $F$ be a divisor on the blow-up $X$ of \pr2 at 
$n \leq 8$ general points, and suppose that 
$ F \cdot  C \geq -1$ for
all exceptional curves $C$ on $X$.  Then 
either $\h{1}{X}{F} = 0$, or $F = r H +  K_X$ with $r \geq 2$, where  $[H]$ 
is primitive and $H$ is smooth, rational, and numerically effective with 
$H^2 = 0$ (in which case $\h{1}{X}{ F} = r-1$).}

\Prf
Since $ F \cdot  C \geq -1$ for all exceptional curves $C$, 
$ C \cdot ( F - K_X) \geq 0$ for all exceptional curves  C, 
which implies $F-K_X = D$ is numerically effective by \ir{eff}.  
Therefore $\h{1}{X}{ F} = \h{1}{X}{D+K_X}= \h{1}{X}{\d{K_X-(D+K_X)}} 
= \h{1}{X}{-D}$, with the center equality due to Serre duality.  
By Ramanujam vanishing \cite{ramanujam} (see \cite{terakawa}
for a characteristic $p$ version, or see Theorem 2.8 of
\cite{birmor}), 
$\h{1}{X}{-D} = 0$ if $D^2 > 0$. If $D^2=0$,
we still have $D\cdot(-K_X)>0$ since, for $n\le8$, $-K_X$ is ample.
Then by an easy calculation applying Lemma 1.4 of \cite{clsrs}
we see that the only possibility is for $[D]$ to be 
in the orbit of $r[L-E_1]$ 
under the action of the Weyl group (see the proof
of \ir{lists}) on $\hbox{Cl}(X)$. I.e., 
$[D]=r[H]$ for some $r$, where $[H]$ is primitive and
$H$ is a smooth rational curve with $H^2 = 0$. 
Now $\h{1}{X}{-D} = r-1$ follows by induction via
$0\to \sh{X}{(-r-1)H}\to \sh{X}{-rH}\to \shn{H}\to 0$.
\qed

We will need to refer to the main result of \cite{algorithm}:

\prclm{Theorem}{Sevenpts}{Let $F$ be a monotone, numerically effective divisor
on the blow up $X$ of \pr2 at 7 general points,
$L,E_1,\cdots,E_7$ being the corresponding exceptional configuration.
Let $t_F$ denote the number of indices $i$ such that
$F\cdot E_i = -F\cdot K_X$ and let $\gamma_F=
\hbox{max}(0,\h{0}{X}{F+L} - 3\h{0}{X}{F})$ (this is the ``expected,''
maximal rank dimension of the cokernel of $\mu_F$).
Then $\dimcok{F}=\hbox{max}(t_F,\gamma_F)$ unless
$F$ is, up to linear equivalence, either
0, $B$, $B-K_X-E_4$, $B-2K_X-E_4-E_5$, 
$B-3K_X-E_4-E_5-E_6$, $B-4K_X-E_4-E_5-E_6-E_7$, $G$  or $G-K_X-E_7$, 
where $B=4L-2E_1-2E_2-2E_3-E_4-\cdots-E_7$ and
$G=5L-2E_1-\cdots-2E_6-E_7$, in which case $\mu_F$ is injective and
$\dimcok{F}=\gamma_F$.}

\Prf
This is Theorem I.6.1 of \cite{algorithm}.
\qed

\prclm{Lemma}{squareZero}{Let 
$C$ be a smooth rational curve on the blow-up $X$ of \pr2 at 8
general points with exceptional configuration
$L,  E_1, \ldots ,  E_8$, and assume $C \cdot L = d$.  Let $0 \leq t$,
and consider the map
$$\mu_{t,C} : \H{0}{X}{\sh{C}{t}} \ot \H{0}{X}{\sh{X}{1}}
\to \H{0}{C}{\sh{C}{t+d}},$$  
given by restriction and multiplication on simple tensors. 
\item{(a)} If $C^2 = -1$ (i.e., $C$ is an exceptional curve on $X$), then
$\mu_{t,C}$ has maximal rank.
\item{(b)} If $C^2 = 0$, then $\mu_{t,C}$ has maximal rank, except
(up to permutation of the indices) in the
following cases: 
\itemitem{(i)} $[C] = [4L - 3E_1 - E_2 - \cdots - E_8]$ and $t=1$, or
\itemitem{(ii)} $[C] = [8L-3E_1-\cdots-3E_7-E_8]$ and $t=3$.
\vskip0in
\noindent For both (i) and (ii), the rank of $\mu_{t,C}$ is one short of 
maximal rank.}

\Prf
It is easy to see that the kernel of the
surjective sheaf map ${\cal O}_C \ot \H{0}{X}{L} \to \sh{C}{d}$  
is $\sh{C}{-a} \oplus \sh{C}{-b}$, for some nonnegative $a$ and 
$b$ with $a+b=d$.  If we assume $a \leq b$,
then it turns out (see the proof of Theorem 3.1 of \cite{curves}) 
that $a \geq \hbox{min}\{d-m_C, m_C\}$, and thus
when $d-m_C$ and $m_C$ differ by at most one, $a$ is the smaller
of $d-m_C$ and $m_C$ and $b$ is the larger. Moreover,
when $a$ and $b$ differ by at most one, then for each $t$ either
$\h{0}{C}{\sh{C}{t-a} \oplus \sh{C}{t-b}}=0$ or
$\h{1}{C}{\sh{C}{t-a} \oplus \sh{C}{t-b}}=0$, hence
$\mu_{t,C}$ has maximal rank for every $t$.

Part (a) now follows by checking \ir{lists}(a)
to see that $d-m_C$ and $m_C$ always 
differ by at most 1 if $C$ is an exceptional curve.

The same proof via \ir{lists}(b) also works for (b), 
except for the following four cases, for which
$d-m_C$ and $m_C$ differ by more than one:
$ 4L - 3E_1 - E_2 - \cdots - E_8$,
$8L - 3E_1 - \cdots - 3E_7 - E_8$,
$10L - 4E_1 - \cdots - 4E_4 - 3E_5 - \cdots - 3E_8$, and
$11L - 4E_1 - \cdots - 4E_7 - 3E_8$.

To handle $C = 4L - 3E_1 - E_2 - \cdots - E_8$, let 
$D = 4L - 3E_1 - E_2 - \cdots - E_7$ (so $D \cdot C = 1$
and $D - C = E_8$), and consider the exact sequence
$$0 \to \sh{X}{D-C} \to \sh{X}{D} \to \sh{C}{1} \to 0.$$
We can compute \dimcok{E_8} trivially
($\mu_{E_8}$ is bijective) and $\mu_D$ has a one-dimensional 
kernel and cokernel by \ir{Sevenpts}, so by 
\ir{mumes}, $\dimker{1,C} = 1$.  Since the
kernel of $\mu_{t,C}$ is $\sh{C}{t-a} \oplus \sh{C}{t-b}$,
we see $$\h{0}{C}{\sh{C}{1-a} \oplus \sh{C}{1-b}} = 1,$$ 
which with $a+b=4$ gives $a=1$ and $b=3$.  
We see that $\mu_{t,C}$ has maximal rank unless
$t=1$, in which case the rank is one short of maximal.

To handle $C = 8L - 3E_1 - \cdots - 3E_7 - E_8$, 
let $D = 8L - 3E_1 - \cdots - 3E_7$ (so again
$D \cdot C = 1$ and $D - C = E_8$).  Consider
the exact sequence
$$0 \to \sh{X}{3D-C} \to \sh{X}{3D} \to \sh{C}{3} \to 0.$$
Note that $3D-C = 2D + E_8$.  
The inclusion $\hbox{ker}(\mu_{2D})\subset \hbox{ker}(\mu_{2D+E_8})$
is clearly an isomorphism since it is induced by the 
canonical injection $\sh{X}{2D}\to \sh{X}{2D+E_8}$
which gives an isomorphism on global sections.
But by \ir{Sevenpts},
$\mu_{2D}$ and $\mu_{3D}$ have kernels of dimensions 0 and 1, 
respectively, so by \ir{mumes}, \dimker{3,C}
is at least 1 dimensional.  Thus $a \leq 3$, but we know $a+b=8$ and
$a$ is at least min$\{d-m,m\}=3$, so in fact $a=3$ and 
$b=5$, which gives the desired result.

To handle $C = 10L - 4E_1 - \cdots - 4E_4 - 3E_5 - \cdots - 3E_8$,
note that $(a,b)$ must be either $(4,6)$ or $(5,5)$.  We just 
need to compute the dimension of the kernel of $\mu_{4,C}$ to 
tell which.  From the exact sequence
$$ 0 \to \sh{X}{-2K_X} \to \sh{X}{-2K_X+C} \to \sh{C}{4} \to 0,$$
and the fact that $\mu_{-2K_X}$ is bijective (see \ir{anexmpl}),
it suffices to compute \dimker{-2K_X+C}.  Let $D = 6L - 3E_1
- 2E_2 - \cdots - 2E_8$ and look at:
$$ 0 \to \sh{X}{-2K_X+C-D} \to \sh{X}{-2K_X+C} \to \sh{D}{2} \to 0.$$
Using the injectivity of $\mu_{2,D}$ from part (a), we reduce
to $\mu_{-2K_X+C-D}$.  Now let $E = 6L - 2E_1 - 3E_2
- 2E_3 - \cdots - 2E_8$, look at 
$$ 0 \to \sh{X}{-2K_X+C-D-E} \to \sh{X}{-2K_X+C-D} \to \sh{E}{2} \to 0,$$
and note $-2K_X+C-D-E = 4L - E_1 - E_2 - 2E_3 - 2E_4 - E_5 - \cdots - E_8$.
Permuting indices gives $F = 4L - 2E_1 - 2E_2 - E_3 - \cdots - E_8$,
which has $l(F) = q(F) = 0$, hence $\mu_{-2K_X+C-D-E}$ and thus
$\mu_{4,C}$ is injective by \ir{qstarlstarlemma}; 
i.e., $a$ must be bigger than 4, so $a=5$.
Again, the desired result follows. 

To handle $C = 11L - 4E_1 - \cdots - 4E_7 - 3E_8$, 
note that we know that $(a,b)$ is either 
$(4,7)$ or $(5,6)$ and as before, the dimension of the kernel of 
$\mu_{4,C}$ tells which.  Let $D=6L - 3E_1
- 2E_2 - \cdots - 2E_8$ and $E=6L - 2E_1 - 3E_2
- 2E_3 - \cdots - 2E_8$ as above, and let $F = 2L - E_3 - \cdots
-E_7$.  Following a process similar to above, we find that
\dimker{-2K_X+C-D-E-F} = \dimker{4,C}, but $-2K_X+C-D-E-F = -K_X$
and we know \dimker{-K_X} = 0 (see \ir{anexmpl}).
Thus $a > 4$ and hence $a=5$, as needed.
\qed

\prclm{Proposition}{partcprop}{Let 
$X$ be a good surface with exceptional configuration
$L, E_1,\ldots,E_8$. Let $F$ be a monotone divisor on $X$
such that $F\cdot E\geq \lambda_E$ for all exceptional
curves $E$ but $F\cdot C< \Lambda_C$ for some exceptional
curve $C$. Then either
\item{(i)} $F \cdot (L - E_1 - E_2)=0$, in which case 
$$\dimcok{F} =
\h{1}{X}{F-(L-E_1)} + \h{1}{X}{F-(L-E_2)},\hbox{ or}$$
\item{(ii)} $\mu_F$ has maximal rank, or
\item{(iii)} $F\cdot E_1=\cdots =F\cdot E_7$ 
(i.e., $F$ is {\rm nearly uniform}), $F$ is numerically effective 
and $ F \cdot  C = 2$, 
where $[C] = [5L - 2E_1 - \cdots - 2E_6 - E_7 - E_8]$.}

\Prf By \ir{eff}, we see that $F$ is effective and
numerically effective, and hence (by (A1)) \h{1}{X}{F} vanishes.
From \ir{lists} and monotonicity of $F$
we can assume that $[C]$ is one of
$[L-E_1-E_2]$, $[3L-2E_1-E_2-\cdots-E_7]$, or $[5L-2E_1-\cdots-2E_6-E_7-E_8]$,
since for all other cases $\lambda_C=\Lambda_C$.
We also may as well assume
that $\mu_F$ fails to have maximal rank.

If $[C]$ is $[L-E_1-E_2]$, then $F\cdot (L-E_1-E_2)=0$ and 
our result follows by \ir{perptodegone}, 
so now we may assume that $F\cdot (L-E_1-E_2)>0$.

If $[C]$ is $[3L-2E_1-E_2\cdots-E_7]$, then $F\cdot C=1$ and, by
\ir{hOneisZero}, either $q^*(F)=h^1(X, F-E_1)=0$, or 
$F=rH+K_X+E_1$ for some $r\geq2$, where
$[H]$ is primitive and $H$ is smooth, rational,
and numerically effective with $C\cdot H=0$ and $H^2=0$. 

In the latter case, 
since $F\cdot (L-E_1-E_2)>0$, $F$ is monotone
and $[H]$ is, up to permutation of the indices, one of the classes
listed in \ir{lists}(b), we see that
$[H]$ can only be $[8L-4E_1-3E_2-\cdots-3E_5-2E_6-2E_7-2E_8]$.
Thus $[F]=[5L-2E_1-\cdots-2E_5-E_6-E_7-E_8]+(r-1)[H]$, which for the
purposes of
induction we denote $F_r$. By a calculation,
$q(F_1)=0$ and $l(F_1)=0$, so $\mu_{F_1}$ 
is injective by \ir{qstarlstarlemma}.
Now consider the exact sequence
$0\to \sh{X}{F_{r-1}}\to \sh{X}{F_r}\to \sh{H}{2}\to 0$.
By \ir{squareZero}, $\mu_{2,H}$ has maximal rank
(and hence here must be injective), so applying
\ir{mumes} and inducting, we see
that $\mu_{F_r}$ is injective for all $r\geq1$,
contradicting our assumption that $\mu_F$ fails to
have maximal rank.

On the other hand, suppose $q^*(F)=0$. We still have
$F\cdot C=1$ and $F\cdot E\geq \lambda_E$ for all
exceptional curves $E$. If $F\cdot E_1>F\cdot E_2$, then
$F-(L-E_1)$ is monotone, and $(F-(L-E_1))\cdot E\ge 0$
for all exceptional curves $E$ unless
$E=5L-2(E_1+\cdots+E_6)-E_7-E_8$ (or one obtained from
this class by permuting the indices)
and $F\cdot E=2$, in which case we at least have
$(F-(L-E_1))\cdot E\ge -1$. As before, by applying \ir{hOneisZero},
either $l^*(F)=0$ or 
$F=(L-E_1)+rH+K_X$ for some $r\ge 2$ with
$[H]$ primitive and $H$ smooth, rational,
and numerically effective with $H^2=0$
If $C\cdot H=0$, then 
$F\cdot C=0$ (contradicting $F\cdot C=1$),
while if $C\cdot H\ge1$, then $F\cdot C\ge r\ge 2$ 
(contradicting $F\cdot C=1$). Thus we see
$l^*(F)=0$ in addition to $q^*(F)=0$, so $\mu_F$ is surjective
by \ir{qstarlstarlemma}, contradicting our 
assumption that $\mu_F$ fails to have
maximal rank.

If, however, $F\cdot E_1=F\cdot E_2$, then 
by checking each possible exceptional curve $E$ we see that either:
\item{(1)} $(F-(L-E_1))\cdot E\ge -1$
for all exceptional curves $E$, or
\item{(2)} $(F-(L-E_1))\cdot E = -2$ for $[E]=[3L-2E_2-E_3-\cdots-E_8]$,
hence $F\cdot E=F\cdot C= 1$ and so 
$F\cdot E_1=\cdots=F\cdot E_8$, or
\item{(3)}  $(F-(L-E_1))\cdot E = -2$ for $[E]=[5L-E_1-2E_2\cdots-2E_7-E_8]$,
hence $F\cdot E=2$ and, since, for $[E']=[5L-2E_1-\cdots-2E_6-E_7-E_8]$,
we have $2=F\cdot E\ge F\cdot E'\ge\lambda_{E'}=2$ by monotonicity, we see
$F\cdot E_1=\cdots=F\cdot E_7$.

In case 1, by applying \ir{hOneisZero},
we see as before that $l^*(F)=0$, since having
$F=(L-E_1)+rH+K_X$ again contradicts $F\cdot C=1$. Thus 
$\mu_F$ is surjective by \ir{qstarlstarlemma},
contrary to assumption.
In case 2, $F=tL-m(E_1+\cdots+E_8)$ for some $t$ and $m$.
Since $F$ is numerically effective, we must have $m\ge 0$ and
$F\cdot (6L-3E_1-2(E_2+\cdots+E_8))\ge 0$, so
$t\ge 17m/6$. Thus $m/2=17m/2-8m \le F\cdot C=1$, so 
$m\le 2$, in which case by \ir{anexmpl} we know 
that $\mu_F$ for the particular $F$ we have here
will have maximal rank, contrary to assumption.
Case 3 is just case (iii) of \ir{partcprop}.

Now we may assume that $F\cdot (3L-2E_1-E_2-\cdots-E_7)>1$,
and we consider the case that $[C]=[5L-2E_1-\cdots-2E_6-E_7-E_8]$
with $F\cdot C=2$. Note now
that $(F-E_1)\cdot E\geq0$ for all exceptional curves
$E$ so $F-E_1$ is, by \ir{eff}, 
numerically effective and effective, 
and hence $q^*(F)=0$.

If $(F-(L-E_1))\cdot E<-1$ for some exceptional curve $E$,
then $E$ can be taken to be $5L-E_1-2E_2-\cdots-2E_6-2E_7-E_8$
and $F\cdot E_1=\cdots=F\cdot E_7$. Otherwise,
we may assume $F\cdot E\geq-1$ for all exceptional curves $E$,
hence by \ir{hOneisZero} either $l^*(F)=0$ (and $\mu_F$ 
is surjective, contrary to assumption), or $F=rH+K_X+L-E_1$ where $r\geq2$,
and $[H]$ is primitive and $H$ is
smooth, rational and numerically effective
with $C\cdot H=0$ and $H^2=0$. 

In the latter case, keeping in mind that $F$
is monotone (which means, since $r\ge 2$, that $H$ must be monotone too),
that $F\cdot C=2$ and that $F\cdot(3L-2E_1-E_2-\cdots-E_7)>1$, 
from \ir{lists} we see that 
$[H]$ must be either $[10L-4E_1-\cdots-4E_4-3E_5-\cdots-3E_8]$
or $[11L-4E_1-\cdots-4E_7-3E_8]$. In each case 
$\mu_F$ ends up having maximal rank. The argument in each
case is similar; here are the details for the latter case.

So let $[H]$ be $[11L-4E_1-\cdots-4E_7-3E_8]$ and 
$[F]=[9L-4E_1-3E_2-\cdots-3E_7-2E_8]+(r-1)[H]$; 
we will denote $(9L-4E_1-3E_2-\cdots-3E_7-2E_8)+tH$ 
by $F_t$, for $t\geq0$.
Consider the exact sequence
$0\to \sh{X}{F_{t-1}}\to \sh{X}{F_t}\to \sh{H}{5}\to 0$.
By \ir{qstarlstarlemma}, $\mu_{F_0}$
is surjective, and by \ir{squareZero}, $\mu_{5,H}$ 
is too. Applying
\ir{mumes} and inducting, we see
that $\mu_{F_t}$ is surjective for all $t\geq0$,
contradicting our assumption that $\mu_F$ fails to have
maximal rank.
\qed

To complete the proof of 
\ir{mainthmgeom}, it still remains to analyze
numerically effective monotone divisors $F$ which are nearly uniform,
have $F\cdot E\ge \lambda_E$ for all exceptional curves $E$ 
and have $F \cdot  C = 2$, 
where $[C] = [5L - 2E_1 - \cdots - 2E_6 - E_7 - E_8]$.
It will be useful to determine all
nearly uniform monotone numerically effective classes;
this is what we do now.
To simplify notation, we will denote the class 
$[d L - a E_1 - \cdots - a E_7 - bE_8]$ by
the triple $(d,a,b)$.

\prclm{Proposition}{conegensprop}{If 
$X$ is the blow-up of \pr2 at eight general points
with exceptional configuration $L,E_1,\ldots,E_8$, 
the classes
$(1,0,0)$, $(3,1,0)$, $(3,1,1)$, $(8,3,0)$, $(8,3,1)$, $(11,4,3)$, 
and $(17,6,6)$ generate the cone of monotone numerically effective 
nearly uniform classes.}

\Prf
Since $(1,0,0)$, $(3,1,0)$, $(3,1,1)$, $(8,3,0)$, $(8,3,1)$,
$(11,4,3)$, and $(17,6,6)$ are all numerically effective, 
monotone and nearly uniform,
any nonnegative \z-linear combination is also numerically
effective, monotone and nearly uniform.

Conversely, let $F=(d,a,b)$ be a nearly uniform class which is
monotone and numerically effective. 
Since $F$ is monotone, we have $a\ge b$, and
since $F$ is numerically effective,
we have $F\cdot E_8\ge 0$, 
$F\cdot (3L-2E_1-E_2-\cdots-E_7)\ge 0$, 
$F\cdot (5L-2E_1-\cdots-2E_6-E_7-E_8)\ge 0$, and
$F\cdot (6L-3E_1-2-\cdots-2E_8)\ge 0$; i.e., we have
$a \ge b$, $b \ge 0$, $3d - 8a \ge 0$, $5d - 13a - b \ge 0$,
and $6d - 15a - 2b \ge 0$.

\noindent It is not hard to check that the rational solution
set to these inequalities is the cone $\Xi({\bf Q})$
given by all nonnegative 
rational linear combinations of $(1,0,0), (8,3,0), (8,3,1), (11,4,3),$
and $(17,6,6)$: each of these classes satisfies all of the inequalities,
but $a = b$ for $(17,6,6)$ and $(1,0,0)$,
$b = 0$ for $(1,0,0)$ and $(8,3,0)$, $3d - 8a =0$ for $(8,3,0)$
and $(8,3,1)$, $5d - 13a - b = 0$ for $(8,3,1)$ and $(11,4,3)$,
and $6d - 15a - 2b = 0$ for $(11,4,3)$ and $(17,6,6)$. Thus we see that
the cone of monotone numerically effective nearly uniform 
classes is just the cone $\Xi=\Xi({\bf Z})$
of integer lattice points in $\Xi({\bf Q})$.

We now show that $\Xi$
is in fact the set of nonnegative 
\z-linear combinations of $(1,0,0)$, $(8,3,0)$, $(8,3,1)$, 
$(11,4,3)$, $(17,6,6)$, $(3,1,0)$ and $(3,1,1)$.  
Let $\langle \ldots\rangle$ denote the cone generated over ${\bf Z}$,
and let $\langle \ldots\rangle_{\bf Q}$ denote the cone 
generated over ${\bf Q}$.
It is easy to see that $\Xi({\bf Q})$ is the union of the rational cones
$$\Xi_1=\langle (11,4,3), (8,3,1), (8,3,0)\rangle_{\bf Q},$$
$$\Xi_2=\langle (11,4,3), (17,6,6), (8,3,0)\rangle_{\bf Q},$$ and
$$\Xi_3=\langle (1,0,0), (17,6,6), (8,3,0)\rangle_{\bf Q}.$$

First, consider $\Xi_1$. Since
$$\left|\matrix{
11 & 8 & 8 \cr
4 & 3 & 3 \cr
3 & 1 & 0\cr}\right| = 1,$$
every integer lattice point which is a rational linear 
combination of $(11,4,3)$, $(8,3,1)$, and $(8,3,0)$ is in fact a 
\z-linear combination; i.e., $\Xi\cap \Xi_1 = 
\langle (11,4,3), (8,3,1), (8,3,0)\rangle$.

Suppose $(d,a,b) \in \Xi\cap\Xi_2$. Then the integer triple $(d,a,b)$ is
$\alpha (17,6,6) + \beta (8,3,0) + \gamma (11,4,3)$ for 
some nonnegative rational numbers $\alpha$, $\beta$ and $\gamma$,
and we have
$$\left[\matrix{
\alpha \cr \beta \cr \gamma \cr}\right]
=
\left[\matrix{
17 & 8 & 11 \cr
6 & 3 & 4 \cr
6 & 0 & 3 \cr}
\right]^{-1}
\left[\matrix{d \cr a \cr b\cr}\right]
=
\left[\matrix{
3 & -8 & -1/3 \cr
2 & -5 & -2/3 \cr
-6 & 16 & 1\cr}\right]
\left[\matrix{d \cr a \cr b \cr}\right].$$
Thus $3d - 8a - b/3  = \alpha$,
$2d - 5a -2b/3 = \beta$, and 
$-6d + 16a + b = \gamma$,
so $\gamma$ is an integer and $\alpha$ and $\beta$ (in lowest terms, of
course) have denominators of 1 or 3.  Thus any element of $\Xi\cap\Xi_2$
is of the form $A+B$, where $A\in \langle (11,4,3), (17,6,6), (8,3,0)\rangle$
and $B=\alpha'(17,6,6)+\beta'(8,3,0)$ with $\alpha',\beta'\in\{0, 1/3, 2/3\}$.
But the only such $B$ which are integer triples are
$(1/3)(17,6,6)+(2/3)(8,3,0)=(11,4,2)=(8,3,1)+(3,1,1)$ and
$(2/3)(17,6,6)+(1/3)(8,3,0)=(14,5,4)=(11,4,3)+(3,1,1)$. This argument shows
that all elements of $\Xi\cap\Xi_2$ are contained in the rational cone
$\langle (11,4,3), (17,6,6), (8,3,0), (8,3,1), (3,1,1) \rangle$.

Finally, suppose $(d,a,b) = \alpha (1,0,0) + \beta (8,3,0) + \gamma (17,6,6)
\in \Xi\cap\Xi_3$; then $6 \gamma \in \z$, $3\beta + 6\gamma \in \z$, and 
$\alpha + 8\beta + 17\gamma \in \z$.  
Thus we may assume $\alpha$ and $\gamma$ are multiples of $1/6$ and
$\beta$ is a multiple of $1/3$.  Thus any element of $\Xi\cap\Xi_3$
is of the form $A+B$, where $A\in \langle (1,0,0), (17,6,6), (8,3,0)\rangle$
and $B=\alpha'(1,0,0)+\beta'(8,3,0)+\gamma'(17,6,6)$ with 
$\alpha',\gamma'\in\{0,1/6,\ldots,5/6\}$ and $\beta'\in\{0, 1/3, 2/3\}$.
By direct check, for every integer triple $B$
we have $B\in \langle (1,0,0), (17,6,6), (8,3,0), (3,1,1), (3,1,0)\rangle$.

Thus $(1,0,0)$, $(3,1,0)$, $(3,1,1)$, $(8,3,0)$, $(8,3,1)$, $(11,4,3)$, 
and $(17,6,6)$ generate the cone of monotone, numerically effective
nearly uniform classes on a blow-up of \pr2\/ at eight general points.
\qed

We now analyze those classes falling into case (iii) of
\ir{partcprop}.

\prclm{Proposition}{specialfamily}{Let 
$F$ be a monotone, numerically effective, nearly uniform
divisor class such that $F\cdot E\ge \lambda_E$ 
for all exceptional curves $E$ and 
$F\cdot(5L - 2E_1 - \cdots - 2E_6 - E_7 - E_8) = 2$. Then either
\item{(a)} $F$ is $(3,1,0) + r(8,3,1)$ 
for some $r \geq 0$, in which case  $\dimcok{F} = r$ and 
$\dimker{F} = r+1$, or
\item{(b)} $\mu_F$ has maximal rank.}

\Prf By \ir{conegensprop}, we know
$$F\in\langle (1,0,0), (3,1,0), (3,1,1), (8,3,0), (8,3,1), (11,4,3), 
(17,6,6)\rangle.$$ Since $F\cdot (5L - 2E_1 - \cdots - 2E_6 - E_7 - E_8)=2$,
one deduces 
that $F$ must be of the form
$H + r(8,3,1) + s(11,4,3)$ for some nonnegative integers $r$ and $s$,
where $H$ is one of
$(3,1,0)$, $(6,2,2)$, $(16,6,0)$, $(20,7,7)$, $(25,9,6)$, or $(34,12,12)$.
(Note that we need not consider the possibility $H = (11,4,1)$ since
$(11,4,1)$ has been accounted for by taking $H = (3,1,0)$ with $r=1$ and 
$s=0$.)  For $F=H+r(8,3,1)+s(11,4,3)$, in order for 
$F\cdot E\ge \lambda_E$ for all exceptional curves $E$ (see
\ir{lists}),
it is easy to check 
that the following additional restrictions
are necessary:
\item{$\bullet$} If $H=(6,2,2)$, then $r>0$, hence we can replace $H=(6,2,2)$ by
$H=(6,2,2)+(8,3,1)=(14,5,3)$ and remove the requirement that $r>0$.
\item{$\bullet$} If $H=(16,6,0)$, then $s>0$, so we replace $H=(16,6,0)$ by
$H=(27,10,3)$.
\item{$\bullet$} If $H=(20,7,7)$, then $r>1$, so we replace $H=(20,7,7)$ by
$H=(36,13,9)$.
\item{$\bullet$} If $H=(34,12,12)$, then $r>2$, so we replace $H=(34,12,12)$ by
$H=(58,21,15)$.

Thus $F$ must be of the form
$H + r(8,3,1) + s(11,4,3)$ for some nonnegative integers $r$ and $s$,
where $H$ is one of
$(3,1,0)$, $(14,5,3)$, $(27,10,3)$, $(36,13,9)$, $(25,9,6)$, or $(58,21,15)$.
We consider each possibility for $H$ in turn, beginning with $H=(3,1,0)$.

So $F = H + r(8,3,1) + s(11,4,3)$, with 
$H = (3,1,0)$. We first consider the case that $s=0$. Note
that $\hbox{cok}(\mu_H)=0$ by \ir{Sevenpts}.
Now, $F\cdot (8,3,1)=3$ so for $r\ge 1$ we have the exact sequence
$0 \to \sh{X}{F-D} \to \sh{X}{F}\to \sh{D}{3} \to 0$ 
with $D = (8,3,1)$. By \ir{squareZero} and an easy calculation, 
$\mu_{3,D}$ has 1-dimensional kernel and cokernel.
Applying \ir{mumes} to the foregoing exact sequence
with $r=1$ (and $s=0$) we see that the induced map
$\hbox{ker}(\mu_F)\to\hbox{ker}(\mu_{3,D})$ is surjective.
Since the restriction of $\sh{X}{D}$ to $D$ is trivial,
the image of the map $\H{0}{X}{H + rD}\to\H{0}{D}{\sh{D}{3}}$
and hence of $\H{0}{X}{H + rD}\otimes\H{0}{X}{L}
\to\H{0}{D}{\sh{D}{3}}\otimes\H{0}{X}{L}$ is the same for all $r\ge1$.
Thus $\hbox{ker}(\mu_F)\to\hbox{ker}(\mu_{3,D})$ is surjective
for all $r\ge1$. We therefore see that the exact sequence
$$\eqalign{0 & \to  \mbox{ker} \mu_{F-D} \to \mbox{ker} \mu_{F} \to \mbox{ker} \mu_{3, D} \cr
             & \to  \mbox{cok} \mu_{F-D} \to \mbox{cok} \mu_{F} \to \mbox{cok} \mu_{3, D} \to 0,\cr}$$
coming from \ir{mumes} is exact separately on kernels
and cokernels. It now follows for $s=0$ and all $r\ge0$
that $\dimker{F} = r+1$ and $\dimcok{F} = r$, as claimed in part
(a).  

Now assume  $F=H + s(11,4,3)$; we find
$F\cdot (11,4,3)=5$ so we have the exact sequence
$0 \to \sh{X}{F-D} \to \sh{X}{F}\to \sh{D}{5} \to 0$ where
this time we take $D = (11,4,3)$.
By \ir{squareZero}, $\mu_{5,D}$ has maximal rank
and one easily checks that $\mu_{5,D}$ therefore is surjective.
Applying \ir{mumes}
and inducting on $s$ we see that $\hbox{cok}(\mu_F)=0$. 

Finally, consider $F=H + s(11,4,3) + r(8,3,1)$ with $s>0$; we find
$F\cdot (8,3,1)=3+s$ so we have the exact sequence
$0 \to \sh{X}{F-D} \to \sh{X}{F}\to \sh{D}{3+s} \to 0$ where
this time we take $D = (8,3,1)$.
By \ir{squareZero}, $\mu_{3+s,D}$ is surjective.
Applying \ir{mumes}
and inducting on $r$ we see that $\hbox{cok}(\mu_F)=0$ for 
all $r\ge0$ and $s>0$. 

Now let $H=(14,5,3)$. Note that $(14,5,3)=(3,1,0)+(11,4,3)$.
Thus $F=H+r(8,3,1)+s(11,4,3)$ is just $(3,1,0)+r(8,3,1)+(s+1)(11,4,3)$,
and our preceding analysis shows that $\hbox{cok}(\mu_F)=0$
for $F=(3,1,0)+r(8,3,1)+(s+1)(11,4,3)$.

The remaining possibilities for $H$ reduce in a similar way
to the case $H=(3,1,0)$ treated above: 
$(27,10,3)$ is $(3,1,0)+3(8,3,1)$; $(36,13,9)$ is 
$(3,1,0)+3(11,4,3)$; $(25,9,6)$ is $(3,1,0)+2(11,4,3)$; and 
$(58,21,15)$ is $(3,1,0)+5(11,4,3)$.
In each instance the reader will easily verify that either (a) or (b) of
the statement
of \ir{specialfamily} is obtained.
\qed

\irrnSection{Examples}{sectionexample}
In this section we show by example how our results give minimal free
resolutions
for fat point subschemes $Z=m_1p_1+\cdots+m_np_n$ of \pr2 with $n\le8$,
where the points
$p_i$ are assumed to be general. We also give two examples for $n>8$
general points on a smooth plane cubic curve, one 
showing that our results sometimes determine resolutions
even though $n>8$ and one showing that sometimes they do not.

If we are interested in a resolution of $I_Z$ for $Z=m_1p_1+\cdots+m_np_n$
with $n<8$ we might as well assume $n=8$ and simply set
some multiplicities $m_i$ equal to 0; i.e., having $n<8$ is no different from
having $n=8$. Now, for our first example, consider $Z=54(p_1+\cdots+p_8)$. 
Recall that $h_Z(t)$ is the Hilbert function of $I_Z$ in degree $t$; thus
$h_Z(t)=\hbox{dim}((I_Z)_t)$. Let $X$ be the 
surface obtained by blowing up the points $p_i$, with the
corresponding exceptional configuration being $L, E_1,\ldots,E_8$. Note that
$D=17L-6(E_1+\cdots+E_8)$ is numerically effective by \ir{eff}.
Denote $tL-54(E_1+\cdots+E_8)$ by $H_t$. Since $H_t\cdot D =
17t-54\cdot6\cdot8$
is negative for $t<153$,
we see $\h{0}{X}{H_t}=0$ for $t<153$.
Since $H_t\cdot E\ge 0$ for all exceptional curves $E$ when $t\ge153$, 
we know $H_t$ is numerically effective for all $t\ge 153$
and hence that $\h{1}{X}{H_t}=0$ and $\h{0}{X}{H_t}=
((H_t)^2-K_X\cdot H_t)/2+1 =
\binom{t+2}{2}-8\binom{54+1}{2}$ for $t\ge 153$.
In other words, $h_Z(t)=0$ for $t<153$, and 
$h_Z(t)=\binom{t+2}{2}-8\binom{54+1}{2}$ for $t\ge153$. 

By vanishing of $h_Z(t)$ for $t<153$
we see $\nu_t=0$ for $t<153$, and,
by vanishing of $\h{1}{X}{H_t}$ for $t\ge153$ and by
the general fact at the end of \ir{anexmpl}
we see $\nu_{t+2}=0$ for all $t\ge 153$.
Since it is obvious that $\nu_{153}=h_Z(153)=55$, we are left with
finding $\nu_{154}$. But $H_{153}\cdot E=0<\lambda_E$ for 
$E = 6L-3E_1-2E_2-\cdots-2E_8$, and likewise for 
$E = 6L-2E_1-3E_2-2E_4-\cdots-2E_8, \ldots$,
so by \ir{mainthmgeom}(b), $\mu_{H_{153}}$ has the same kernel
as the maps corresponding to $H_{153}-(6L-3E_1-2E_2-\cdots-2E_8)$, 
$H_{153}-(6L-3E_1-2E_2-\cdots-2E_8)-(6L-2E_1-3E_2-2E_4-\cdots-2E_8)$, etc., 
and we find eventually that $\mu_{H_{153}}$ has the same kernel
as $\mu_H$ with $H=9L-3E_1-\cdots-3E_8$. By \ir{mainthmgeom}(a),
$\mu_H$ has maximal rank, and since $\h{0}{X}{H}=7$ and $\h{0}{X}{H+L}=18$
by \ir{eff},
we see that $\mu_H$ must be surjective with 3-dimensional kernel.
Thus $\mu_{H_{153}}$ has 3-dimensional kernel, and from 
$\h{0}{X}{H_{153}}=55$ and $\h{0}{X}{H_{154}}=210$, we see
$\nu_{154}$ equals 48 (rather than the maximal rank value of 45). 

Using the relation $\nu_t-s_t=\Delta^3h_Z(t)$, we find that
$s_t$ is 0 for $t<154$ or $t>155$, $s_{154}=3$ and $s_{155}=99$.
Thus the minimal free resolution of $I_Z$ is
$0\to F_1\to F_0\to I_Z\to 0$ where $F_0=R^{55}[-153]\oplus R^{48}[-154]$
and $F_1=R^{3}[-154]\oplus R^{99}[-155]$.

Consider now another example. Suppose $Z=mp_1+\cdots+mp_n$ for $n\ge 9$
general points of a smooth cubic curve in \pr2. We show our results 
recover the resolution of $I_Z$, which is known in this case 
(see section 3.2.1 of \cite{resolution}). So let $H_t=tL-m(E_1+\cdots+E_n)$
with $m>0$ and let $D_X$ as usual be a smooth section of 
$-K_X$, where $X$ is obtained from \pr2 by blowing up the points $p_i$.
If $t<3m$, then $H_t\cdot(-K_X)<0$, so $\h{0}{X}{H_t}=\h{0}{X}{H_t+K_X}$,
but we still have $(H_t+K_X)\cdot (-K_X)<0$ and iterating we eventually
find that $\h{0}{X}{H_t}=\h{0}{X}{H_t+mK_X}=0$ (the last equality follows
since $L\cdot (H_t+mK_X)<0$). Thus $h_Z(t)=0$ for $t<3m$, hence $\mu_t(Z)$
has maximal rank for $t<3m$. For $t=3m$, then $|H_t|=\{mD_X\}$, so 
$\mu_t(Z)$ has maximal rank. For $t>3m$, then $H_t=(t-3m)L-mK_X$, so 
$E\cdot H_t = (t-3m)E\cdot L - mE\cdot K_X > E\cdot L > \Lambda_E$
for every exceptional curve $E$ and either $H_t\cdot(-K_X)\ge2=\Lambda_{D_X}$
(hence $\mu_t(Z)$ has maximal rank by \ir{mainthmgeom}(a))
or $H_t\cdot(-K_X)<2=\lambda_{D_X}$ (hence $\mu_{H_t}$ and  
$\mu_{H_t+K_X}$ have kernels of the same dimension 
by \ir{mainthmgeom}(b), and iterating, for some $l$ we eventually
obtain $H_t+lK_X$ falling under case (a)
of \ir{mainthmgeom}). Thus in any case we
can compute the rank of $\mu_{H_t}$ for 
every $t$, which makes it easy to work out the resolution for any
particular $m$ and $n$.

Finally, consider $Z=156p_1+121(p_2+\cdots+p_7)+104p_8+78p_9$.
Again, let $H_t=tL-(156E_1+108(E_2+\cdots+E_7)+104E_8+78E_9))$.
It turns out that $[27L-12E_1-9(E_2+\cdots+E_7)-8E_8-6E_9]$ 
is the class $[E]$ of an exceptional curve $E$, and that 
$[H_t]=[(t-351)L+13E]$. From this it is easy to check that
$\mu_{H_t}$ is injective for $t\le 351$ (since $\h{0}{X}{H_t}=0$
for $t<351$ while $\h{0}{X}{H_{351}}=1$) and (by \ir{anexmpl},
since $\h{1}{X}{H_t}=0$ for $t>351$) surjective for $t>352$.
However, for $t=352$ we have $H_t\cdot C\ge \Lambda_C$ for all $C\in\Gamma_X$
except $C=E$, for which we have $\lambda_E=12<H_t\cdot E = 14 < 15 =
\Lambda_E$,
and hence \ir{mainthmgeom} does not apply and, indeed, the rank
of $\mu_{H_t}$ is not known.

\References

\bibitem{\catalisano} M. V. Catalisano, {\it Linear Systems of Plane Curves 
through Fixed ``Fat'' Points of \pr2}, 
J.\ Alg.\  142 (1991), 81-100.

\bibitem{\curves} S. Fitchett,
{\it On Bounding the Number of Generators for Fat Point
Ideals on the Projective Plane}, to appear, Journal of Algebra.

\bibitem{\sfmaps} S. Fitchett, {\it Maps of linear systems
on blow ups of the projective plane}, 
to appear, J. Pure and Applied Alg.

\bibitem{\clsrs} B. Harbourne, {\it Complete Linear Systems
on Rational Surfaces}, Trans. Amer. Math. Soc., 289(1985), 213--226.

\bibitem{\duke} B. Harbourne, {\it Blowings-up of \pr2
and their blowings-down}, Duke.\ Math. J., 52(1985), 129--148.

\bibitem{\mtnwest} B. Harbourne, {\it Rational surfaces with $K^2>0$}, 
Proc. Amer. Math. Soc. 124 (1996), 727--733.

\bibitem{\ars} B. Harbourne, {\it Anticanonical 
rational surfaces}, Trans. 
Amer. Math. Soc. 349 (1997), 1191--1208.

\bibitem{\birmor} B. Harbourne, {\it Birational morphisms of
rational surfaces}, J. Alg. 190 (1997), 145--162.

\bibitem{\algorithm} B. Harbourne, {\it An Algorithm for Fat Points on \pr2}, 
Canad. J. Math. 52 (2000), 123--140.

\bibitem{\resolution} B. Harbourne, {\it Free Resolutions of Fat Point 
Ideals on \pr2}, J. Pure and Applied Alg. 125 
(1998), 213--234.

\bibitem{\igp} B. Harbourne, {\it The Ideal Generation 
Problem for Fat Points}, J. Pure and Applied Alg. 145 
(2000), 165--182.

\bibitem{\survey} B. Harbourne, 
{\it Problems and Progress:  A survey on fat points in \pr2},
to appear, Queen's Papers in Pure and Applied Math.
(http://www.math.unl/$\sim$bharbour/Survey.tex).

\bibitem{\mumford} D. Mumford,
{\it Varieties defined by quadratic equations,}
Questions on algebraic varieties,
Corso C.I.M.E. 1969, Rome: Cremoneses, 1970, 30-100.

\bibitem{\nagata} M. Nagata, {\it On rational surfaces, II}, 
Mem.\ Coll.\ Sci.\ 
Univ.\ Kyoto, Ser.\ A Math.\ 33 (1960), 271--293.

\bibitem{\ramanujam} C. P. Ramanujam, {\it Supplement
to the article ``Remarks on the Kodaira
vanishing theorem''}, J. Indian Math. Soc. {\bf 38} (1974), 121--124.

\bibitem{\terakawa} H. Terakawa, {\it The $d$-very ampleness
on a projective surface in characteristic $p$}, 
Pac. J. Math. 187 (1999), 187--199.

\bye